\documentclass[12pt]{article}
\usepackage{amsmath}
\usepackage{amsfonts}
\usepackage{amssymb}

\usepackage{amsthm,amsmath,amssymb}
\usepackage{verbatim}
\usepackage{array}
\usepackage[colorlinks=true,citecolor=black,linkcolor=black,urlcolor=blue]{hyperref}

\newtheorem{theorem}{Theorem}
\newtheorem{lemma}[theorem]{Lemma}

\newtheorem{claim}[theorem]{Claim}

\newtheorem{conjecture}[theorem]{Conjecture}

\newtheorem{remark}[theorem]{Remark}
\theoremstyle{remark}

\usepackage{lipsum}

\begin{document}
\title{{\Large On Ramsey numbers of $3$-uniform Berge cycles}}
\author{\small   Leila Maherani$^{\textrm{b},1}$ \ and \  Maryam Shahsiah$^{\textrm{a},\textrm{b},2}$\\
		\small  $^{\textrm{a}}$ Department of Mathematics, Khansar Campus, University of Isfahan, Isfahan, Iran\\	
	\small  $^{\textrm{b}}$School of Mathematics, Institute for
	Research
	in Fundamental Sciences (IPM),\\
	\small  P.O.Box: 19395-5746, Tehran,
	Iran\\	
	\small \texttt{E-mails: 
		maherani@ipm.ir, m.shahsiah@khn.ui.ac.ir}}
\date {}
\maketitle \footnotetext[1] {\tt This research was supported by a grant from IPM.} \vspace*{-0.5cm}
\footnotetext[2] {\tt This research was in part supported by a grant from IPM (No. 99050423).} \vspace*{-0.5cm}

%
\setlength{\parindent}{0pt}
\maketitle
\begin{abstract}
For an arbitrary graph $G$, a hypergraph $\mathcal{H}$  is called Berge-$G$ if there is a bijection $\Phi :E(G)\longrightarrow E(
\mathcal{H})$ such that for each $e\in E(G)$, we have $e\subseteq \Phi (e)$.  We
denote  by $\mathcal{B}^rG$, the family of $r$-uniform Berge-$G$ hypergraphs. 
For families  $\mathcal{H}_1, \mathcal{H}_2,\ldots, \mathcal{H}_t$  of $r$-uniform hypergraphs,  the  Ramsey number  $R(\mathcal{H}_1, \mathcal{H}_2,\ldots, \mathcal{H}_t)$   is
the smallest integer $n$ such that in every $t$-hyperedge coloring of  $\mathcal{K}_{n}^r$ there is 
a monochromatic copy of a hypergraph in $\mathcal{H}_i$ of color $i$, for some $1\leq i\leq t$.
 Recently, the Ramsey problems of Berge hypergraphs have been studied by many  researchers.\\
 In this paper, we focus on  Ramsey number involving $3$-uniform Berge cycles and we prove that for $n \geq 4$, $ R(\mathcal{B}^3C_n,\mathcal{B}^3C_n,\mathcal{B}^3C_3)=n+1.$  Moreover, for  $m  \geq n\geq 6$ and $m\geq 11$, we show that $R(\mathcal{B}^3K_m,\mathcal{B}^3C_n)= m+\lfloor \frac{n-1}{2}\rfloor -1.$ 
This is the first result of Ramsey number for two different families of Berge hypergraphs.\\
\end{abstract}
\textbf{Keywoeds}: Extremal combinatorics, Ramsey theory, Hypergraphs.
\section{Introduction}

A \textit{hypergraph} $\mathcal{H}$ is a pair $\mathcal{H}= (V, E)$, where $V$ is a finite and non-empty set of vertices
and $E$ is a collection of distinct non-empty subsets of $V$ (the set of hyperedges).
 The set of vertices and the set of hyperedges of $\mathcal{H}$ are denoted by  $V(\mathcal{H})$ and $E(\mathcal{H})$, respectively.
An {\it $r$-uniform
	hypergraph} is a hypergraph such that all of its hyperedges have size $r$. Every graph is a $2$-uniform hypergraph. For a vertex $v\in V(\mathcal{H})$ the \textit{degree} of $v$ denoted by $d_{\mathcal{H}}(v)$ and is the number of hyperedges in $E(\mathcal{H})$ containing $v$. Let $S\subseteq V(\mathcal{H})$. The set 
  $N_{S}(v):=\{u\in S: \exists\ h \in E(\mathcal{H}), \{v,u\}\subseteq h\}$ stands for the neighborhood of $v$ in $S$. 
  We also denote by $\mathcal{H}[S]$ the subhypergraph of $\mathcal{H}$ induced on $S$.
   By the $\mathcal{H}-v$, we mean  the subhypergraph of $\mathcal{H}$ made by removing $v$ and all the hyperedges containing it.
  A \textit{complete $r$-uniform hypergraph} of order $n$, denoted by $\mathcal{K}_n^r$ , is a hypergraph consisting of all the $r$-subsets of a set $V$ of cardinality $n$.

For given $n \geq 3$ and $r \geq 2$, a
\textit{complete $r$-uniform Berge hypergraph} of order $n$, denoted by $\mathcal{K}_n^{(r)}$, is an $r$-uniform hypergraph
with the main vertices $\{v_1, v_2, \ldots , v_n\}$  and ${n \choose 2}$ 
distinct hyperedges $h_{ij}$, $1 \leq i < j \leq n$, each of size $r$,
where every $h_{ij}$ contains both $v_i$ and $v_j$. An \textit{$r$-uniform Berge cycle} of length $n$,  denoted by $\mathcal{C}_n^{(r)}$, is an
$r$-uniform hypergraph with the core sequence $v_1, v_2, \ldots , v_n$,  as the main vertices, and distinct
hyperedges $h_1, h_2, \ldots , h_n$ such that $h_i$ contains  $v_i$ and $v_{i+1}$, where addition is done modulo
$n$.\\
For an arbitrary graph $G$, the definition of Berge-$G$ was introduced by Gerbner and Palmer \cite{{gerbdef}}.
 A hypergraph $\mathcal{H}$  is called \textit{Berge-$G$} if there is a bijection $\Phi :E(G)\longrightarrow E(
\mathcal{H})$ such that for each $e\in E(G)$, we have $e\subseteq \Phi (e)$. 
In other words, one can
enlarge edges of $G$ to obtain $\mathcal{H}$, or  can shrink  hyperedges of $\mathcal{H}$ to obtain a copy of $G$.
A graph $G$ may have several non-isomorphic Berge copies and $\mathcal{H}$ may be Berge hypergraph of several different graphs.  We
denote  by $\mathcal{B}^rG$, the family of $r$-uniform Berge-$G$ hypergraphs. Clearly,  $\mathcal{K}_n^{(r)}$ and  $\mathcal{C}_n^{(r)}$ are Berge-$K_n$ and Berge-$C_n$, respectively. Therefore,  for every $n\geq 3$ and $r \geq 2$, $\mathcal{K}_n^{(r)} \in \mathcal{B}^rK_n$ and  $\mathcal{C}_n^{(r)} \in \mathcal{B}^rC_n$, where $K_n$ is a complete graph of order $n$ and $C_n$  is a cycle  of length $n$.\\
 
 By a \textit{$t$-hyperedge coloring} of $\mathcal{H}$, we mean a function $C:E(\mathcal{H})\longrightarrow \{\mathfrak{c}_1,\mathfrak{c}_2,\ldots,\mathfrak{c}_t\}$. 
 For $1\leq j\leq t$,  let $\mathcal{H}_{\mathfrak{c}_j}$ be the  subhypergraph of $\mathcal{H}$ induced by the hyperedges of color $\mathfrak{c}_j$.
 For a hypergraph $\mathcal{H}$, \textit{shadow graph} $\Gamma(\mathcal{H})$  is a graph  with $V(\Gamma(\mathcal{H}))=V(\mathcal{H})$ and two vertices are adjacent  if and only if they are incident to  a  hyperedge of $\mathcal{H}$. For a $t$-hyperedge coloring of $\mathcal{H}$,  the graph $\Gamma_i(\mathcal{H})$ is defined such a way that $V(\Gamma _i(\mathcal{H}))=V(\Gamma(\mathcal{H}))$, $E(\Gamma _i(\mathcal{H}))=E(\Gamma(\mathcal{H}))$ and  for every edge $e\in E(\Gamma_i(\mathcal{H}))$, $c(e)$ is a list of colors such that $\mathfrak{c}_j \in c(e)$, $1\leq j\leq t$, if $e$ is contained in at least $i$ hyperedges of  $\mathcal{H}_{\mathfrak{c}_j}$.\\
 
  Let  $\mathcal{H}_1, \mathcal{H}_2,\ldots, \mathcal{H}_t$ be families of $r$-uniform hypergraphs. The  \textit{Ramsey number}  $R(\mathcal{H}_1, \mathcal{H}_2,$\\
 $\ldots, \mathcal{H}_t)$   is
 the smallest integer $n$ such that in every $t$-hyperedge coloring of  $\mathcal{K}_{n}^r$ there is 
 a monochromatic copy of a hypergraph in $\mathcal{H}_i$ of color $\mathfrak{c}_i$, for some $1\leq i\leq t$.
 We use the notation  $R(\mathcal{H};t)$ if $\mathcal{H}_i\cong \mathcal{H}$, $1\leq i\leq t$.  For a family $\mathcal{F}$ of $r$-uniform hypergraphs, we say that the hypergraph $\mathcal{H}$ is $\mathcal{F}$-free if $\mathcal{H}$   contains no member of $\mathcal{F}$ as a subhypergraph.\\

\begin{remark}
 Throughout the paper, to avoid the confusion, we will always apply edge for  graph edge  and hyperedge for hypergraph edge. Also, for the simplicity, we denote by letters ${\bf r}$, ${\bf b}$ and ${\bf g}$ the colors \text{red}, \text{blue} and \text{green}, respectively. Therefore, in a $2$-hyperedge coloring of the hypergraph $\mathcal{H}$ with colors red and blue, each edge $e \in E(\Gamma _i(\mathcal{H}))$ has a list in form $\{\textbf{r}\}$, $\{\textbf{b}\}$ or $\{\textbf{r},\textbf{b}\}$.  Also,  in a $3$-hyperedge coloring of the hypergraph $\mathcal{H}$ with colors red, blue and green, each edge $e \in E(\Gamma _i(\mathcal{H}))$ has a list in form $\{\textbf{r}\}$, $\{\textbf{b}\}$, $\{\textbf{g}\}$, $\{\textbf{r},\textbf{b}\}$, $\{\textbf{r},\textbf{g}\}$, $\{\textbf{b},\textbf{g}\}$ or $\{\textbf{r},\textbf{b},\textbf{g}\}$.
 Moreover, we denote by $\mathcal{H}_{\textbf{r}}$,  $\mathcal{H}_{\textbf{b}}$ and $\mathcal{H}_{\textbf{g}}$ the  subhypergraphs induced by the hyperedges of colors red, blue and green, respectively.
\end{remark} 
 

 For a graph $G$ and two vertices $u,v \in V(G)$, we use the notation $u\sim _G v$ if $u$ and $v$ are adjacent vertices in $G$.
 A \textit{matching} in a simple graph $G$ is a set of  edges with no shared endpoints.
  The vertices incident to the edges of a matching $M$
  are saturated by $M$. For $X\subseteq V(G)$, by a matching $M$ of $X$ we mean that $X$ is saturated by $M$. Let $\mathcal{A} = \{A_1,\ldots, A_n\} $ be a collection of subsets of a set $X$. A \textit{system of
  distinct representatives} (briefly, SDR) for $A$ is a set of distinct elements $a_1, \ldots, a_n$ in $X$ such
  that $a_i \in A_i$, $1\leq i\leq n$.\\


%
 

 The extremal problems of Berge hypergraphs is one of the important problems in combinatorics which has been of interest to many investigators.
  Recently,  Tur\'{a}n number of some families of Berge hypergraphs for instance  $\mathcal{B}^r C_n$, $\mathcal{B}^r P_n$, $\mathcal{B}^r K_n$, $\mathcal{B}^r T_n$, $\mathcal{B}^r K_{r,s}$, ect. has received considerable attention. To see more results  we refer the reader to \cite{gyori, katona, davoodi, mahcomp, gyarfcomp, gerbcomp, palmer2, vizer}.\\
  
 The study of  Ramsey numbers  of Berge hypergraphs was initiated by Gy\'{a}rf\'{a}s et al. \cite{gyarfconj1}, in 2008.  They conjectured that $R(\mathcal{B}^r C_n; r-1) = n$ 
 for large enough $n$ and proved it for $r = 3$ and $n \geq 5$.
 Gy\'{a}rf\'{a}s, S\'{a}rk\"{o}zy and Szem\'{e}r\'{e}di \cite{zemeredi}   used the method of regularity lemma to prove the asymptotic form of this conjecture. More precisely, they   showed that  $R(\mathcal{B}^r C_n;r-1)  = \textbf{(}1 + o(1)\textbf{)}n$.
 In \cite{longcycle},  the authors   showed that for $n \geq 140$, in every $3$-hyperedge coloring of $\mathcal{K}^4_n$
 there
 is a monochromatic Berge-cycle of length at least $n - 10$.  Regarding the latter case, Maherani and Omidi  \cite{mahcycle} gave an affirmative answer to this conjecture for $r=4$ by proving  $R(\mathcal{B}^4C_n;3) = n $ for $n\geq 75$. More recently,  this conjecture has been proved by  Omidi \cite{omidicycle}.
 
 Using the method of regularity lemma, Gy\'{a}rf\'{a}s and S\'{a}rk\"{o}zy established the following theorem on $3$-color Ramsey number of Berge cycles.
 \begin{theorem}\em{\cite{gyarfconj2}}\label{gyarf}
 For all $\eta >0$ there exists $n_0$ such that for every $n>n_0$ the following holds.
 $$R(\mathcal{B}^3C_n;3)\leq (\dfrac{5}{4}+\eta)n.$$
 \end{theorem}
 They also posed the following conjecture. 
 \begin{conjecture}\em{\cite{gyarfconj2}}\label{gyarf2}
For all $\eta>0$ and positive integer $r$, there exists $n_0=n_0(\eta,r)$ such that for every $n>n_0$, every $r$-hyperedge coloring of $\mathcal{K}_n^r$ contains a monochromatic Berge cycle of length at least $(\dfrac{2r-2}{2r-1}-\eta)n$.
 \end{conjecture}

In this paper, we consider the $3$-color Ramsey number of $3$-uniform Berge cycles and prove the following theorem.
\begin{theorem}\label{ccc}
	For every $n\geq 4$,  $ R(\mathcal{B}^3C_n,\mathcal{B}^3C_n,\mathcal{B}^3C_3)=n+1.$
\end{theorem}

 Recently, The  Ramsey number of $\mathcal{B}^3K_n$, has been investigated by many researchers \cite{axenovich,omidi,salia}.  Axenovich and  Gy\'{a}rf\'{a}s \cite{axenovich} studied $R(\mathcal{B}^3 K_n;t)$  for $n \in  \{3, 4\}$.   Gerbner et al. \cite{omidi} claimed that  $R
  (\mathcal{B}^3 K_n;t)$ is bounded by a
  polynomial of $n$ depending on $t$. 
   In \cite{omidi},  the authors also determined  the exact value of the Ramsey number of $3$-uniform Berge trees.
     Furthermore, they found a general  bounds for the Ramsey number of $3$-uniform Berge hypergraphs.
   Salia et al. \cite{salia} classified the exact values of $R(\mathcal{B}^3 K_m,\mathcal{B}^3 K_n)$ based on the  parameters $m$ and $n$. Also, for large $n$, a superlinear lower bound was obtained for $R(\mathcal{B}^3 K_n;3)$. It is worth mentioning that the Ramsey number  $R(\mathcal{B}^r K_n;t)$ has been studied as well, for the higher uniformity.  For more details, see \cite{salia, gerbram, omidi, palvolgy}. \\
   
   As mentioned above, the Ramsey numbers  of Berge cycles and also  complete Berge hypergraphs have been studied by many researchers. For the first time, in this paper,   the combination of these two families  has been considered.  More precisely, we prove the following theorem.
   \begin{theorem}\label{ck}
   		Let  $n\in \{4,5\}$ and $m\geq 8$ be integers.  
   	 Then $R(\mathcal{B}^3K_m,\mathcal{B}^3C_n)=m+1$. Moreover, for  $m  \geq n\geq 6$ and $m\geq 11$, we have $R(\mathcal{B}^3K_m,\mathcal{B}^3C_n)=m+\lfloor \frac{n-1}{2}\rfloor -1$.
   \end{theorem}

For the reader's convenience, we summarise here how the paper is organised.
In section 2, we will prove Theorem \ref{ccc}. This section includes two subsections. In the first one, Theorem \ref{ccc} will be affirmed for small cases $n\in\{4,5\}$ (Theorem \ref{c4c4c3} and \ref{c3c5c5}). In the second one, using induction, we will prove Theorem \ref{ccc} for general case $n \geq 6$. We close the paper by establishing the proof of  Theorem \ref{ck}  in section 3. Theorem \ref{ck} will be proven for small cases $n\in\{4,5\}$ (Theorem \ref{ck-small})  and then for the general case (Theorem \ref{ck-general}). 

\section{The 3-color Ramsey number of $\mathcal{B}^3C_n$}

Axenovich and Gy\'{a}rf\'{a}s \cite{axenovich} determined the exact value of  $R(\mathcal{B}^3C_3;t)$ for $t \in \{2,3,4,5,6,8\}$. In this section,  we focus on   3-color Ramsey number of Berge cycles where lengths of two cycles may be infinitely grown. 
 Using the   Tur\'{a}n number of $\mathcal{B}^3C_3$ and   some applicable lemmas, we first  prove Theorem \ref{ccc} for  small cases $n \in \{4,5\}$ and then apply induction to prove it in  general case. 
In the following, we show that for  $n\geq 4$,  $ R(\mathcal{B}^3C_n,\mathcal{B}^3C_n,\mathcal{B}^3C_3)\geq n+1$. To see that, consider two sets $A$ and $B$, where $|A|=2$ and $|B|=n-2$. Now, color  all triples intersecting $A$ in exactly two points by green and  the remaining triples by red and blue, arbitrarily. Clearly, this coloring contains no red or blue  $\mathcal{C}_n^{(3)}$ and no green $\mathcal{C}_3^{(3)}$. Hence, 
  \begin{equation}\label{lowerbound}
  R(\mathcal{B}^3C_n,\mathcal{B}^3C_n,\mathcal{B}^3C_3)\geq n+1.
  \end{equation}
\subsection{Small cases}
 The following upper bound for the Tur\'{a}n number of $\mathcal{B}^3C_3$ has been found by Gy\H{o}ri \cite{gyorik3}.  He also showed that this bound is tight.
\begin{theorem}\em{\cite{gyorik3}}\label{exc3}
	For $n \geq 3$, $ex(n, \mathcal{B}^3C_3)\leq  \dfrac{n^2}{8}$.
\end{theorem}

In \cite{axenovich}, it was shown that $R(\mathcal{B}^3C_3;3)=5$.
In the sequel, we demonstrate that Theorem \ref{ccc} holds for $n \in \{4,5\}$. First, we prove  the case $n=4$.  
\begin{theorem}\label{c4c4c3}
  $ R(\mathcal{B}^3C_4,\mathcal{B}^3C_4,\mathcal{B}^3C_3)=5$.	
\end{theorem}
 \noindent\textbf{proof. }By the relation (\ref{lowerbound}), it is sufficient to show that    $ R(\mathcal{B}^3C_4,\mathcal{B}^3C_4,\mathcal{B}^3C_3)\leq 5$.
	Consider a $3$-hyperedge coloring of  $\mathcal{K}_5^3$ with colors red, blue and green. Clearly, at least four hyperedges, say  $h_1, h_2, h_3 , h_4$ receive the same color.  If they are  green, using Theorem \ref{exc3}, the assertion holds. Otherwise,  we may assume that these hyperedges are  red and $|h_1\cap h_2|=2$. By a simple case analysis, it can be proven that no matter how  are $h_3$ and $h_4$, we have a red $\mathcal{C}_4^{(3)}.$  
	$\hfill\blacksquare$\\
	
 According to the proof  of Theorem \ref{c4c4c3}, the following statement holds.
 \begin{remark}\label{k5-6e}
 	 Every $3$-uniform hypergraph with five vertices  and four hyperedges contains a $\mathcal{C}_4^{(3)}$.
 \end{remark}
Before we prove  the case $n=5$ of Theorem \ref{ccc} , we present the following  lemma.
 
\begin{lemma}\label{k5-3e}
	 Every $3$-uniform hypergraph with five vertices  and seven hyperedges contains a $\mathcal{C}_5^{(3)}$ as a subhypergraph.
\end{lemma}
 \noindent\textbf{proof. }Suppose that $\mathcal{K}=\mathcal{K}_5^3$ with $V(\mathcal{K})=\{v_1,v_2,v_3,v_4,v_5\}$ and $E(\mathcal{K})=\{h_1,h_2,\ldots, h_{10}\}$. Let $\mathcal{H}=\mathcal{K}\setminus \{h_1,h_2,h_3\}.$   Without loss of generality, we may assume that $h_1=\{v_1,v_4,v_5\} $ and  $h_2=\{v_2,v_4,v_5\} $. If $h_3=\{v_3,v_4,v_5\} $, then $v_1,v_2,v_5,v_3,v_4$ is the core sequence of a $\mathcal{C}_5^{(3)}$ with the hyperedges  $\{v_1,v_2,v_4\},\{v_2,v_3,v_5\},\{v_1,v_3,v_5\},\{v_2,v_3,v_4\},\{v_1,v_3,v_4\}$ in  $\mathcal{H}$.
	Otherwise, if $h_3=\{v_2,v_3,v_4\}$ or  $h_3=\{v_1,v_2,v_3\}$, then $v_1,v_2,v_3,v_4,v_5$ represents the core sequence of a $\mathcal{C}_5^{(3)}$ with  hyperedges  $\{v_1,v_2,v_5\},\{v_2,v_3,v_5\},\{v_1,v_3,v_4\},\{v_3,v_4,v_5\},\{v_1,v_3,v_5\}$ in $\mathcal{H}$.
Therefore, we may assume that  $h_3=\{v_1,v_2,v_5\}$. In this case, $\mathcal{H}$ contains  a $\mathcal{C}_5^{(3)}$ with the core sequence $v_1,v_2,v_3,v_4,v_5$ and the  hyperedges  $\{v_1,v_2,v_3\},\{v_2,v_3,v_4\},\{v_1,v_3,v_4\},$\\$\{v_3,v_4,v_5\},\{v_1,v_3,v_5\}$, so we are done.
	$\hfill\blacksquare$\\
\begin{theorem}\label{c3c5c5}
	$ R(\mathcal{B}^3C_5,\mathcal{B}^3C_5,\mathcal{B}^3C_3)=6$.
	\end{theorem}
 \noindent\textbf{proof. }The lower bound holds by the relation (\ref{lowerbound}). Assume that $\mathcal{H}=\mathcal{K}_6^3$ is  $3$-hyperedge colored with colors red, blue and green such that $ \mathcal{H}_\textbf{g}$ is   $\mathcal{B}^3C_3$-free. By Theorem \ref{exc3},  $\mathcal{H}_\textbf{g}$ contains at most four hyperedges. Therefore, $|E (\mathcal{H}_\textbf{r}) |\geq 8$ or $|E (\mathcal{H}_\textbf{b}) |\geq 8$. Without loss of generality,  suppose that $|E(\mathcal{H}_\textbf{r}) |\geq 8$. Let $V(\mathcal{H})=\{v_1,v_2,v_3,v_4,v_5, v_6\}$. We show that there is a set $S \subset V(\mathcal{H})$ such that  $|S|=5$ and $\mathcal{H}_\textbf{r}[S]$ has  at least four  hyperedges. Suppose to the contrary that there is no such a set $S$. Counting the set $T=\{(h,E)\ : \  h\in E(\mathcal{H}_\textbf{r}),\ h\subset E,\ |E|=5\}$ in two sides implies that $|T|\leq 18$ and $| T| \geq 24$, a contradiction. Let $S=\{v_1,v_2,v_3,v_4,v_5\}$ be such a set. Using Remark \ref{k5-6e}, $\mathcal{H}_\textbf{r}[S]$ contains  a  $\mathcal{C}_4^{(3)}$ as a subhypergraph, say $\mathcal{C}$. Assume that $v_1,v_2,v_3,v_4$ and $h_1,h_2,h_3,h_4$ are the core sequence and the  hyperedges  of  $\mathcal{C}$, respectively. Note that   $d_{\mathcal{H}_\textbf{r}}(v_6)\geq 2$.	 Otherwise, $\mathcal{H}_\textbf{r}[S]$ contains seven hyperedges and by Lemma \ref{k5-3e}, the proof is completed.
 By considering the following cases, we have a $\mathcal{C}_5^{(3)}$ in $\mathcal{H}_{\textbf{r}}$.
\begin{itemize}
	\item[(i)] There are at least three  hyperedges in $\mathcal{H}_{\textbf{r}}$ containing both  $v_5$ and $v_6$.
	
	Clearly, there are two consecutive vertices $v_i$ and $v_{i+1}$ on $\mathcal{C}$  for which two hyperedges  $\{v_i,v_5,v_6\}$ and $\{v_{i+1}, v_5,v_6\}$ are red. Hence, Berge cycle $\mathcal{C}$ can be extended to  a  $\mathcal{C}_5^{(3)}$ in $\mathcal{H}_\textbf{r}$ by embedding $v_6$ between   $v_i$ and $v_{i+1}$ on $\mathcal{C}$, removing $h_i$ and adding$\{v_i,v_5,v_6\}$ and $\{v_{i+1}, v_5,v_6\}$.  
	
	\item[(ii)] There are exactly two hyperedges in $\mathcal{H}_{\textbf{r}}$ containing both    $v_5$ and $v_6$.
	
	If there are two consecutive vertices $v_i$ and $v_{i+1}$ on $\mathcal{C}$ such that hyperedges $\{v_i,v_5,v_6\} $ and $\{v_{i+1},v_5,v_6\}$ are red, by an argument  similar to the previous case, we are done. Therefore,  we may assume that   the hyperedges $\{v_1,v_5,v_6\} $ and $\{v_3,v_5,v_6\}$ are red. As mentioned, we have $|E(\mathcal{H}_\textbf{r}) |\geq 8$. On the other hand, the set $\{v_1,v_2,v_3,v_4\}$ contains at most four red hyperedges. This implies that there are at least two red hyperedges including exactly one of the $v_5$ or $v_6$. Without loss of generality, suppose that there is at least one red hyperedge including $v_5$ and not containg $v_6$. 
	If $\{v_1,v_3,v_5\}$ is red, then the sequence $v_1,v_4,v_3,v_6,v_5$ with the hyperedges $h_4,h_3,\{v_3,v_5,v_6\},\{v_1,v_5,v_6\},\{v_1,v_3,v_5\}$ make a  $\mathcal{C}_5^{(3)}$ in $\mathcal{H}_\textbf{r}$. Otherwise, there is a red hyperedge in form $\{v_j,v_{j'},v_5\}$, such that  $v_j \in \{v_2,v_4\}$.
	If $v_j=v_2$, then $v_1,v_2,v_5,v_3,v_4$ represents the core sequence of a red $\mathcal{C}_5^{(3)}$ with the hyperedges $h_1, \{v_2,v_5,v_{j'}\}, \{v_3,v_5,v_6\}, h_3,h_4$. Otherwise $v_j=v_4$. In this case, $v_1,v_2,v_3,v_5,v_4$ with the hyperedges $h_1, h_2, \{v_3,v_5,v_6\}, \{v_4,v_5,v_{j'}\}, h_4$ make a  red $\mathcal{C}_5^{(3)}$.
	 
	\item[(iii)] There is only one  hyperedge  in $\mathcal{H}_{\textbf{r}}$ containing both   $v_5$ and $v_6$.
	
	The hyperedge $q=\{v_5,v_6,v_t\}$ is red for exactly one vertex $v_t \in \{v_1,v_2,v_3,v_4\}$. Since $d_{\mathcal{H}_\textbf{r}}(v_6)\geq 2$,	  there exist $v_j,v_{j'} \in \{v_1,\ldots,v_4\}$ such that the hyperedge $\{v_j,v_{j'},v_6\}$ is red. If there are two hyperedges $g_1,g_2 \in E(\mathcal{H}_\textbf{r})$ and two consecutive vertices $v_i, v_{i+1} \in \{v_1,v_2,v_3,v_4\}$  where $\{v_i,v_6\}\subset g_1$ and $\{v_{i+1},v_6\}\subset g_2$,  then  Berge cycle $\mathcal{C}$ can be extended to  a  $\mathcal{C}_5^{(3)}$ in $\mathcal{H}_\textbf{r}$ by embedding $v_6$ between   $v_i$ and $v_{i+1}$ on $\mathcal{C}$, removing $h_i$ and adding $g_1$ and $g_2$.  
Therefore, $q$ and $\{v_6,v_t,v_{t+2}\}$ are the only  red hyperedges containing $v_6$. On the other hand, since $|E (\mathcal{H}_\textbf{r}) |\geq 8$, there are two red hyperedges containing $v_5$ and not containing $v_6$.\\
	  Let $f_1=\{v_5,v_{t-1},v_t\}$,  $f_2=\{v_5,v_t,v_{t+1}\}$ and $f_3=\{v_5,v_{t-1},v_{t+1}\}$. If  $f_1$ is red, then we can extend $\mathcal{C}$ to  a red $\mathcal{C}_5^{(3)}$   by embedding $v_5$ between $v_{t-1}$ and $v_t$  on $\mathcal{C}$, removing $h_{t-1}$  and adding   hyperedges   $f_1$ and $q$. Similarly, if  one of the $f_2$ or $f_3$ is red, we have a  $\mathcal{C}_5^{(3)}$ in $\mathcal{H}_\textbf{r}$   by embedding $v_5$ between $v_{t}$ and $v_{t+1}$. 
	   Therefore, we may suppose that  $f_1$, $f_2$ and $f_3$ are blue and at least two hyperedges of $f_4=\{v_5,v_{t-1},v_{t-2}\}$, $f_5=\{v_5,v_t,v_{t+2}\}$ or $f_6=\{v_5,v_{t+1},v_{t+2}\}$ are red.\\
	     If  $f_4$ and $f_5$  are red, then by setting $v_5$ between $v_{t-1}$ and $v_{t+2}$, removing $h_{t-2}$ and using $f_4$ and $f_5$   we have a  red $\mathcal{C}_5^{(3)}$. Also, if $f_5$ and $f_6$ are red,  we can make a    red $\mathcal{C}_5^{(3)}$ by embedding $v_5$ between $v_{t+1}$ and $v_{t+2}$, removing $h_{t+1}$ and using $f_5$ and $f_6$.  Otherwise, the hyperedges $f_4$ and $f_6$ are red and $f_5$ is blue. Thus,  $\mathcal{H}_{\textbf{r}}[\{v_1,v_2,v_3,v_4\}]\cong \mathcal{K}_4^3$ and we may assume that $\mathcal{C}\cong \mathcal{K}_4^3$. Obviously,  one can easily embed $v_5$ between $v_{t+1}$ and $v_{t+2}$ to obtain a   $\mathcal{C}_5^{(3)}$ in $\mathcal{H}_{\textbf{r}}$.
	 
	\item[(iv)] There is no  hyperedge in $\mathcal{H}_{\textbf{r}}$ containing both  $v_5$ and $v_6$.
	
	Clearly, there are two consecutive vertices $v_i$ and $v_{i+1}$ on $\mathcal{C}$  such that two hyperedges in form $\{v_i, v_{\ell},v_6\}$ and $\{v_{i+1}, v_{\ell'},v_6\}$  are red for some  $ \ell, \ell' \in \{1,\ldots,4\}$. Hence, Berge cycle $\mathcal{C}$ can be extended to  a  $\mathcal{C}_5^{(3)}$ in $\mathcal{H}_{\textbf{r}}$ by embedding $v_6$ between   $v_i$ and $v_{i+1}$ on $\mathcal{C}$, removing $h_i$ and adding $\{ v_i, v_{\ell},v_6\}$ and $\{ v_{i+1}, v_{\ell'},v_6\}$.  
	$\hfill\blacksquare$\\
\end{itemize}
	
\subsection{Proof of Theorem \ref{ccc}}

 Before we present the proof of Theorem \ref{ccc}, we will prove the following lemma. 

\begin{lemma}\label{exten}
		Consider a $t$-hyperedge coloring of   $3$-uniform hypergraph $\mathcal{H}$  with colors $\{\mathfrak{c}_1,\ldots,\mathfrak{c}_t\}$.
	If $C_n \subseteq \Gamma_2(\mathcal{H})$ and  there is a color $\mathfrak{c}_{\ell}$ for some  $1\leq \ell \leq t$ such that for every $e \in E(C_n)$,  $\mathfrak{c}_{\ell}\in c(e)$, then $\mathcal{C}_n^{(3)} \subseteq  \mathcal{H}_{\mathfrak{c}_{\ell}}$.	
\end{lemma}
\noindent\textbf{proof. }Suppose that $v_1v_2\ldots v_n$ is a cyclic order of vertices of $C_n$.
	 Let  $\mathcal{A}_i$ be the set
	of hyperedges of color $\mathfrak{c}_{\ell}$ in $\mathcal{H}$  containing $v_i$ and  $ v_{i+1}$, $1\leq i \leq n$, ($v_{n+1}=v_1$). Now, let $G^*=[X,Y]$ be the bipartite graph defined as follows. Let
	 $X=\{v_iv_{i+1} \ :\ 1\leq i \leq n \}$, $Y=\cup_{i=1}^n \mathcal{A}_i$, and for each $1\leq i\leq n$ and $A \in Y$, $v_iv_{i+1}\sim_{G^*}A$,  if and only if $A$ contains  $v_i$ and $v_{i+1}$. Since $C_n \subseteq \Gamma_2(\mathcal{H})$,   $d_{G^*}(v_iv_{i+1}) \geq 2$ for every $1\leq i \leq n$. Also,   $d_{G^*}(A) \leq 2$ for every $A \in Y$. Using the Hall's theorem,  $G^*$ contains a matching of $X$ that leads to  a  $\mathcal{C}_n^{(3)}$ of color  $\mathfrak{c}_{\ell}$ in $\mathcal{H}$.
$\hfill\blacksquare$\\

{\bf Proof of Theorem \ref{ccc}}:
By the relation (\ref{lowerbound}), it  suffices to prove 
 $R(\mathcal{B} ^3C_n,\mathcal{B}^3C_n,\mathcal{B}^3C_3)\leq n+1$. For this purpose, we use   induction on $n$.
	 For the  cases $n \in {\{4,5\}}$, it is true by Theorems \ref{c4c4c3} and \ref{c3c5c5}. Suppose that $n \geq 6$ and  $R(\mathcal{B} ^3C_{n'},\mathcal{B}^3C_{n'},\mathcal{B}^3C_3)=n'+1$ for $4\leq n' <n$. 
Let  $\mathcal{H}=\mathcal{K}_{n+1}^3$  be $3$-hyperedge colored with colors red, blue and green.  Assume that $ \mathcal{H}_\textbf{g}$ is   $\mathcal{B}^3C_3$-free.  We show that $\mathcal{C}_n^{(3)} \subseteq \mathcal{H}_\textbf{r}$ or   $\mathcal{C}_n^{(3)} \subseteq \mathcal{H}_\textbf{b}$.
%

  First, we prove the following two claims.
  
 \begin{claim}\label{claim1-c3cmcm}
For every pair of edges $e=uv$ and $e'=uv'$ in $\Gamma_2(\mathcal{H})$, none of the following cases hold.
\begin{itemize}
\item [\rm{(i)}] If ${\bf g}\in c(e)\cap c(e')$, then there is a vertex $w \in V(\mathcal{H})\setminus \{u,v,v'\}$ such that both of hyperedges $\{w,u,v\}$ and $\{w,u,v'\}$ are green.
	\item[\rm{(ii)}]  $c(e)=\{{\bf g}\}$ and ${\bf g} \in c(e').$
\end{itemize}
 \end{claim} 
 	\noindent\textbf{proof of Claim \ref{claim1-c3cmcm}. }(i) Suppose not. So, there is a vertex $w \in V(\mathcal{H})\setminus \{u,v,v'\}$ such that both of hyperedges $h=\{w,u,v\}$ and $h'=\{w,u,v'\}$ are green. Since  ${\bf g}\in c(e)\cap c(e')$, there is a hyperedge $h''\neq h$  in $\mathcal{H}_\textbf{g}$ such that $ \{u,v\} \subset h''$. Therefore, we have a green $\mathcal{C}_3^{(3)}$ with core sequence $u, w, v$ and hyperedges $h', h,h''$, a contradiction. \\
 	
 	(ii) Suppose to the contrary that  $c(e)=\{{\bf g}\}$ and ${\bf g} \in c(e')$. Since $c(e)=\{{\bf g}\}$, for at most two vertices $x$ and $x'$ hyperedges $\{u,v,x\}$ and $\{u,v,x'\}$ are not green. The hypothesis $n \geq 6$ implies that there are at least three hyperedges in $\mathcal{H}_\textbf{g}$ containing $u$ and $v$.
 	  If $h=\{u,v,v'\} $ is green, then there are distinct green hyperedges $h'$ and  $h''$, such that  $h'\neq h $, $h''\neq h$,  $\{u,v'\}\subset h'$ and $\{u,v\}\subset h''$. Hence, $u,v,v'$ represents the core sequence of a  $\mathcal{C}_{3}^{(3)}$ with the hyperedges $h'',h,h'$, a contradiction.
 	   Otherwise,  $h $ is red or blue. In this case, by the assumption $n\geq 6$, one can see that there is a vertex $w$ such that hyperedges $\{w,u,v\}$ and $\{w,u,v'\}$ are green. This contradicts the  Case (i).
  	$\hfill\square$\\

Let $V(\mathcal{H})=\{v_1,v_2, \ldots , v_{n+1}\}$. 
By the induction hypothesis, removing every vertex from $\mathcal{H}$ leads to have a red or blue $\mathcal{C}_{n-1}^{(3)}$.  Without loss of  generality, let $V'=\{v'_1,v'_2, \ldots , v'_{\lceil\frac{n+1}{2}\rceil}\}$ be the set of vertices for which  the  subhypergraph $\mathcal{H}-v'_{i}$  contains a red  $ \mathcal{C}_{n-1}^{(3)}$,   $1\leq i \leq \lceil\frac{n+1}{2}\rceil$. Also, suppose that $v_{n+1}\in V'$ and  $\mathcal{C}= \mathcal{C}_{n-1}^{(3)}$ is the  Berge cycle  with the core sequence  $v_1,v_2, \ldots, v_{n-1}$  and the hyperedges $h_1,h_2, \ldots , h_{n-1}$ in $\mathcal{H}_{\textbf{r}}-v_{n+1}$. If there is a red $\mathcal{C}_{n}^{(3)}$ in $\mathcal{H}$, we are done. Therefore, we may assume that $ \mathcal{H}_\textbf{r}$ is   $\mathcal{B}^3C_n$-free.

\begin{claim}\label{eiej}
	The  following two statements hold.
	\begin{itemize}
		\item [\rm{(i)}] For every edge $e_i=v_{n+1}v_i$, $1\leq i\leq n$, in the graph $\Gamma_2(\mathcal{H})$,  $c(e_i)\neq \{{\bf r}\}$.
		\item [\rm{(ii)}] Let $e_i=v_{n+1}v_i$ and $e_j=v_{n+1}v_j$, $1\leq i\neq j \leq n$, be edges of $\Gamma_2(\mathcal{H})$. If  $c(e_i)=\{{\bf r},{\bf g}\}$, then $c(e_j)\neq \{{\bf r},{\bf g}\}$.
	\end{itemize}
\end{claim}

\noindent\textbf{proof of Claim \ref{eiej}. }(i) Suppose to the contrary that $c(e_i)= \{{\bf r}\}$. First, let $1\leq i\leq n-1$.  Assume that the  hyperedge $f_1=\{v_{n+1},v_{i-1},v_{i}\}$  is red. Since there is a red hyperedge $f_2\neq f_1$ contaning $v_{n+1}$ and $v_i$,   we have a  $\mathcal{C}_{n}^{(3)}\subseteq \mathcal{H}_\textbf{r}$, by setting $v_{n+1}$ between  $v_{i-1}$ and $v_{i}$, removing hyperedge $h_{i-1}$  and adding  hyperedges $f_1$ and  $f_2$, a contradiction. Therefore, we may assume that $\{v_{n+1},v_{i-1},v_{i}\}$ is blue and $\{v_{n+1},v_i,v_{i+1}\}$ is green.\\
	Since $n \geq 6$ and $c(e_i)= \{{\bf r}\}$, there exist two consecutive vertices $v_t$ and $v_{t+1}$ on $\mathcal{C}$, $t \notin \{i-2,i-1,i,i+1\}$ such that $\{v_{n+1},v_i,v_{t}\}$ and $\{v_{n+1},v_i,v_{t+1}\}$ are red. Clearly we have a red  $\mathcal{C}_{n}^{(3)}$, by embedding $v_{n+1}$ between $v_t$ and $v_{t+1}$, removing $h_t$ and using hyperedges $\{v_{n+1},v_i,v_{t}\}$ and $\{v_{n+1},v_i,v_{t+1}\}$, that is a contradiction. 
	Now, let  $i=n$. Since $c(e_n)=\{\textbf{r}\}$, there are at most two vertices $v_j$ and $v_{j'}$ on $\mathcal{C}$ such that hyperedges $\{v_{n+1},v_n,v_{j}\}$ and $\{v_{n+1},v_n,v_{j'}\}$ are of colors blue and green. The fact $n\geq 6$ implies that there are two consecutive vertices $v_t$ and $v_{t+1}$ on $\mathcal{C}$  where $\{v_{n+1},v_n,v_{t}\}$ and $\{v_{n+1},v_n,v_{t+1}\}$ are red. Hence, one can  find a  $\mathcal{C}_{n}^{(3)}$ in $\mathcal{H}_\textbf{r}$   by an argument  similar to the above, a contradiction.\\
		
	(ii) 
	Suppose to the contrary that $c(e_j)= \{\textbf{r},\textbf{g}\}$.
	First, assume that $1\leq i, j\leq n-1$. 
	We may assume that $v_i$ and $v_j$ are not consecutive vertices on $\mathcal{C}$. Otherwise, let  $v_i=v_t$ and $v_j=v_{t+1}$. Since there are two distinct  hyperedges $f_1$ and $f_2$ in $\mathcal{H}_\textbf{r}$ such that $\{v_{n+1},v_t\}\subset f_1$ and $\{v_{n+1},v_{t+1}\}\subset f_2$,  $\mathcal{C}$ can be extended to a red   $\mathcal{C}_{n}^{(3)}$ by embedding $v_{n+1}$ between $v_t$ and $v_{t+1}$, removing $h_t$ and adding $f_1$ and $f_2$, a contradiction. \\
	 None of the  hyperedges $\{v_{n+1},v_i,v_{i-1}\}$, $\{v_{n+1},v_i,v_{i+1}\}$, $\{v_{n+1},v_j,v_{j-1}\}$ and $\{v_{n+1},v_j,v_{j+1}\}$  are  red,  by an argument similar to the Case (i). Hence,  since $\textbf{b}\notin c(e_i)$ and $\textbf{b}\notin c(e_j)$,
	one of the hyperedges $\{v_{n+1},v_i,v_{i-1}\}$  or $\{v_{n+1},v_i,v_{i+1}\}$ and also one of the hyperedges $\{v_{n+1},v_j,v_{j-1}\}$  or $\{v_{n+1},v_j,v_{j+1}\}$ are green.
	Without loss of generality, we may assume that $\{v_{n+1},v_i,v_{i-1}\}$ and  $\{v_{n+1},v_j,v_{j+1}\}$ are green ( the argument for the other cases is similar).\\
	  Moreover,  none of the hyperedges $\{v_{n+1},v_j,v_{i-1}\}$, $\{v_{n+1},v_j,v_{i+1}\}$,   $\{v_{n+1},v_i,v_{j-1}\}$  and $\{v_{n+1},v_i,v_{j+1}\}$ are red. Otherwise, suppose that the hyperedge $\{v_{n+1},v_j,v_{i-1}\}$  is red. Since there is a red hyperedge containing $v_{n+1}$ and $v_i$,  we can embed $v_{n+1}$ between $v_{i-1}$ and $v_i$  to construct a red $\mathcal{C}_{n}^{(3)}$, a contradiction. \\
	  By the assumption $\textbf{b}\notin c(e_j)$, at least one of the $\{v_{n+1},v_j,v_{i-1}\}$  or $\{v_{n+1},v_j,v_{i+1}\}$ is green.
	If $\{v_{n+1},v_i,v_{i+1}\}$ is blue, then $\{v_{n+1},v_i,v_{j+1}\}$ is green. This is a contradiction to the Case (i) of Claim \ref{claim1-c3cmcm} by setting $w=v_{j+1}$. Otherwise,  $\{v_{n+1},v_i,v_{i+1}\}$ is green. Since $\{v_{n+1},v_j,v_{i-1}\}$  or $\{v_{n+1},v_j,v_{i+1}\}$ is green, we can find a vertex $w\in \{v_{i-1},v_{i+1}\}$ such that hyperedges  $\{v_{n+1},v_i,w\}$ and $\{v_{n+1},v_j,w\}$ are green. This contradicts  the Case (i) of the Claim \ref{claim1-c3cmcm}, as well.
	
	Now, let $j=n$. By an argument similar to the above,  the hyperedges  $\{v_{n+1},v_i,v_{i-1}\}$, $\{v_{n+1},v_i,v_{i+1}\}$, $\{v_{n+1},v_n,v_{i-1}\}$ and $\{v_{n+1},v_n,v_{i+1}\}$  are not red. On the other hand, by the Case (i) of the Claim  \ref{claim1-c3cmcm}, there is no vertex $w$ on  $\mathcal{C}$ such that both of the  hyperedges $\{v_{n+1},v_i,w\}$ and $\{v_{n+1},v_n,w\}$ are green.
	Therefore, We may assume that $\{v_{n+1},v_i,v_{i-1}\}$ and $\{v_{n+1},v_n,v_{i+1}\}$ are blue and $\{v_{n+1},v_i,v_{i+1}\}$ and $\{v_{n+1},v_n,v_{i-1}\}$ are green. 
	Since  $\textbf{b}\notin c(e_i)$ and $\textbf{b}\notin c(e_n)$, all hyperedges in form $\{v_{n+1},v_i,v_{\ell}\}$ and $\{v_{n+1},v_n,v_{\ell}\}$, $1\leq \ell \leq n-1$, $\ell \notin \{i-1,i+1\}$, are red or green.
	 If one of the  hyperedges $\{v_{n+1},v_i,v_{i-2}\}$ or $\{v_{n+1},v_n,v_{i-2}\}$ is red, then  $\{v_{n+1},v_i,v_{i-3}\}$ and $\{v_{n+1},v_n,v_{i-3}\}$ are green. This contradicts the Case (i) of the Claim  \ref{claim1-c3cmcm} by setting $w=v_{i-3}$. Otherwise,  both of the hyperedges $\{v_{n+1},v_i,v_{i-2}\}$ and $\{v_{n+1},v_n,v_{i-2}\}$ are green. It is also a contradiction to the Case (i) of the Claim  \ref{claim1-c3cmcm} by setting $w=v_{i-2}$.
	$\hfill\square$\\

Let $E'$ be the set of edges of $\Gamma_2(\mathcal{H})$ such that for every $e\in E'$, $\textbf{b}\in c(e)$. Also, let $\Gamma'(\mathcal{H})$ be the spaning  subgraph of $\Gamma_2(\mathcal{H})$ induced by $E'$.
By the Claims \ref{claim1-c3cmcm} and \ref{eiej}, there are at least $n-1$ vertices $\{v_{i_1}, v_{i_2}, \ldots, v_{i_{n-1}}\}\subset \{v_1,v_2,\ldots, v_n\}$ such that  $v_{n+1}v_{i_j}\in E(\Gamma'(\mathcal{H}))$ for each $1\leq j\leq n-1$. In other words, $d_{\Gamma'(\mathcal{H})}(v_{n+1})\geq n-1$. It is reminded that  $V'=\{v'_1,v'_2, \ldots , v'_{\lceil\frac{n+1}{2}\rceil}\}$ is the set of vertices for which   subhypergraph $\mathcal{H}_{\textbf{r}}-v'_{i}$  contains a    $\mathcal{C}_{n-1}^{(3)}$,   $1\leq i \leq \lceil\frac{n+1}{2}\rceil$. Therefore, 
similarly, we have $d_{\Gamma'(\mathcal{H})}(v'_i)\geq n-1$, for every $v'_i \in V'$, $1 \leq i\leq  \lceil\frac{n+1}{2}\rceil$. 
Set  $V''= V(\mathcal{H})\setminus V'$. We may assume that  $ V''=\{v''_1,v''_2, \ldots , v''_{\lfloor{\frac{n+1}{2}} \rfloor}\}$.
In the following, we prove that the graph  $\Gamma'(\mathcal{H})$ contains a cycle $C_n$. Then, using Lemma \ref{exten}, we have     $\mathcal{C}_{n}^{(3)} \subseteq  \mathcal{H}_\textbf{b}$, that completes the proof.

First, let $n$ be odd. Clearly,  $|V'|=|V''|=\frac{n+1}{2}$. Assume that  $v'_1,v'_2, \ldots, v'_{ \frac{n+1}{2}}$ are ordered vertices such that $v'_1\sim _{\Gamma'(\mathcal{H})}v'_{ \frac{n+1}{2}} $ and $|N_{V''}(v'_1)\cap N_{V''}(v'_2)|\geq |N_{V''}(v'_{\ell})\cap N_{V''}(v'_{\ell'})|$, for $1\leq \ell \neq \ell' \leq \frac{n+1	}{2}$, $\{\ell,\ell'\} \neq \{1,2\}$.
We show that there are  distinct vertices $x_1,x_2, \ldots , x_{\frac{n+1}{2} -1} $ in $ V''$  such that for every $i$, $v'_i \sim_{\Gamma'(\mathcal{H})}x_i \sim_{\Gamma'(\mathcal{H})}v'_{i+1}$.  Clearly, the cyclic ordered vertices  $$v'_1,x_1,v'_2,x_2,v'_3,\ldots, v'_{\frac{n+1}{2}-1}, x_{\frac{n+1}{2}-1},v'_{\frac{n+1}{2}} $$ 

represents a cycle $C_n$ in $\Gamma'(\mathcal{H})$. 
  For $1\leq i \leq \frac{n+1}{2}-1$, set $$A_{i,i+1}= \left \{v''_j\ :\ v'_i \sim_{\Gamma'(\mathcal{H})}v''_j \sim_{\Gamma'(\mathcal{H})}v'_{i+1}, 1\leq j \leq \frac{n+1}{2} \right \},$$
  and $$\mathcal{A}=\left \{A_{1,2}, A_{2,3}, \ldots, A_{\frac{n+1}{2}-1,\frac{n+1}{2}}\right \}.$$ 
   Since for every $i$   vertex $v'_i \in V'$ is not adjacent to  at most one  vertex of $(V' \cup V'')\setminus \{v'_i\}$ in $\Gamma'(\mathcal{H})$,  $|A_{i,i+1}|\geq |V''|-2 \geq \frac{n+1}{2}-2$.
  Every  SDR for the family  $\mathcal{A} $ indicates desired vertices $x_1,x_2, \ldots , x_{\frac{n+1}{2} -1}.$ In the sequel, we show that $\mathcal{A}$ contains an SDR.\\
   Let $T=\{A_{i,i+1} \in \mathcal{A} \ :\ |A_{i,i+1}|= \frac{n+1}{2}-2\}$ and $T'=\mathcal{A}\setminus T$.
  If $|T|=0$, then for every $A_{i,i+1} \in \mathcal{A}$, $|A_{i,i+1}|\geq \frac{n+1}{2}-1$. Since $|\mathcal{A}| = \frac{n+1}{2}-1$, it is easy to see that the Hall's theorem guarantees the existence of  an SDR for $\mathcal{A}$. 
  If $|T'|=0$, then for every $A_{i,i+1} \in \mathcal{A}$, $|A_{i,i+1}|= \frac{n+1}{2}-2$. By the properties of ordering of $v'_1,v'_2, \ldots, v'_{ \frac{n+1}{2}}$, we conclude that  $|N_{V''}(v'_{\ell})\cap N_{V''}(v'_{\ell'})| = \frac{n+1}{2}-2$, for every $1\leq \ell,\ell' \leq \frac{n+1	}{2}$. Then, we may assume that $\{v''_i,v''_{i+1}\}\nsubseteq A_{i,i+1}$, $1\leq i \leq \frac{n+1	}{2}-1$. Hence, $v''_3,v''_4,\ldots,v''_{\frac{n+1}{2}},v''_1$ is an SDR for $\mathcal{A}$.
  Otherwise, we may assume that $|T|=\ell $, $1\leq \ell \leq \frac{n+1}{2}-2$,  and $|T'|=\ell' = \frac{n+1}{2}-1-\ell$. Therefore, we can assume that $T=\{{A}'_{1,2},\ldots, {A}'_{\ell,\ell+1}\}$ and $T'=\{{A}''_{1,2},\ldots, {A}''_{\ell',\ell'+1}\}.$ First, we show that $T$ has an SDR. For every $1\leq t\leq \ell$,
  	 $$ \left|{A}'_{i_1,i_1+1}\cup \ldots\cup {A}'_{i_t,i_t+1}\right|\geq \frac{n+1}{2}-2 \geq \ell\geq t.$$
  	  Using the Hall's theorem, $T$ contains an SDR, say $x'_1,x'_2, \ldots , x'_{\ell}$. Set $X'=\{x'_1,x'_2, \ldots , x'_{\ell}\}$. Now, let
  	 $${B}_{i,i+1}= {A}''_{i,i+1}\setminus X'.\ \ \ \ 1\leq i\leq \ell' $$
  	 For every $1\leq t'\leq \ell'$,
  	  $$ \left|{B}_{i_1,i_1+1}\cup \ldots\cup {B}_{i_{t'},i_{t'}+1}\right|\geq \frac{n+1}{2}-1- |X'|=\frac{n+1}{2}-1- \ell=\ell'\geq t'.$$
  	  Hence, $T'$ contains an SDR, say $x''_1,x''_2, \ldots , x''_{\ell'}$. Clearly, the set $ \{x'_1,x'_2, \ldots , x'_{\ell},x''_1,x''_2, \ldots , x''_{\ell'}\} $ consists an SDR for $\mathcal{A}$.\\ 

Now, let $n$ be even. Then $|V'|=\lceil\frac{n+1}{2}\rceil$ and $|V''|=\lceil\frac{n+1}{2}\rceil-1$.
Assume that vertices  $v'_1,v'_2, \ldots, v'_{ \lceil\frac{n+1}{2}\rceil}$ are ordered  such that $v'_{ 1}\sim _{\Gamma'(\mathcal{H})}  v'_{ \lceil\frac{n+1}{2}\rceil}\sim _{\Gamma'(\mathcal{H})}v'_{\lceil\frac{n+1}{2}\rceil-1} $ and   $|N_{V''}(v'_1)\cap N_{V''}(v'_2)|\geq |N_{V''}(v'_{\ell})\cap N_{V''}(v'_{\ell'})|$, for $1\leq \ell,\ell' \leq \lceil\frac{n+1	}{2}\rceil$, $\{\ell,\ell'\} \neq \{1,2\}$.
We show that there are  distinct vertices $x_1,x_2, \ldots , x_{\lceil\frac{n+1}{2}\rceil -2}$ in  $ V''$  such that for every $i$, $v'_i \sim_{\Gamma'(\mathcal{H})}x_i \sim_{\Gamma'(\mathcal{H})}v'_{i+1}$. 
Similar to the   previous case, for $1\leq i \leq \lceil\frac{n+1}{2}\rceil-2$,
 set 
 $$A_{i,i+1}=\left\lbrace v''_j\ :v'_i \sim_{\Gamma'(\mathcal{H})}v''_j \sim_{\Gamma'(\mathcal{H})}v'_{i+1}, 1\leq j \leq \lceil\frac{n+1}{2}\rceil -1 \right\rbrace ,$$
and $$\mathcal{A}=\left\lbrace A_{1,2}, A_{2,3},\ldots, A_{\lceil\frac{n+1}{2}\rceil-2,\lceil\frac{n+1}{2}\rceil-1} \right\rbrace. $$ 
Since for every $i$,  vertex $v'_i \in V'$ is not adjacent to at most one  vertex of $(V' \cup V'')\setminus \{v'_i\}$ in $\Gamma'(\mathcal{H})$,  we have  $|A_{i,i+1}|\geq |V''|-2\geq  \lceil\frac{n+1}{2}\rceil-3$.
Every  SDR for family  $\mathcal{A} $ indicates desired vertices $x_1,x_2, \ldots , x_{\lceil\frac{n+1}{2}\rceil -2}.$\\
  Let $T=\{A_{i,i+1} \in \mathcal{A} \ :\ |A_{i,i+1}|= \lceil\frac{n+1}{2}\rceil-3\}$ and $T'=\mathcal{A}\setminus T$. Since $|V'|> |V''|$ and by the properties of ordering $v'_1,\ldots, v'_{ \lceil\frac{n+1}{2}\rceil}$,  $|N_{V''}(v'_1)\cap N_{V''}(v'_2)|\geq \lceil\frac{n+1}{2}\rceil -2$. Hence,   $|T'|\neq 0$.
   If $|T|=0$, then for every $A_{i,i+1} \in \mathcal{A}$, $|A_{i,i+1}|\geq \lceil\frac{n+1}{2}\rceil-2$. Since $|\mathcal{A}| = \lceil\frac{n+1}{2}\rceil-2$, one can easily see that the Hall's theorem guarantees the existence of  an SDR for $\mathcal{A}$. Now, let  $|T|=\ell $, $1\leq \ell \leq \lceil\frac{n+1}{2}\rceil-3$,  and $|T'|=\ell' = \lceil\frac{n+1}{2}\rceil-2-\ell$.
  By an argument similar to the previous case, there is an SDR $x_1,x_2,\ldots, x_{\lceil\frac{n+1}{2}\rceil-1} $ for $\mathcal{A}$.
  The cyclic ordered vertices  
  $$v'_1,x_1,v'_2,x_2,v'_3,\ldots, v'_{\lceil\frac{n+1}{2}\rceil-2}, x_{\lceil\frac{n+1}{2}\rceil-2},v'_{\lceil\frac{n+1}{2}\rceil-1}, v'_{\lceil\frac{n+1}{2}\rceil} $$
   represents a cycle $C_n$ in $\Gamma'(\mathcal{H})$  and the proof is completed.
$\hfill\blacksquare$\\
 \section{Ramsey number of $\mathcal{B}^3K_m$ versus  $\mathcal{B}^3C_n$  }
  In this section, the Ramsey number    $ R(\mathcal{B}^3K_m, \mathcal{B}^3C_n) $ is established for every $m\geq n\geq 3$ and $m\geq 11$. For this purpose, we first consider the small values of $n\in \{3,4,5\}$ and $m\geq n$ (Theorem \ref{ck-small}) and then determine  $ R(\mathcal{B}^3K_m, \mathcal{B}^3C_n)$ for   $m\geq n\geq 6$ and $m\geq 11$  (Theorem \ref{ck-general}). \\

 Given a family $\mathcal{F}$ of $r$-uniform hypergraphs, the {\it Tur\'{a}n number} of $\mathcal{F}$ for a given
 positive integer $N$, denoted by $ex(N,\mathcal{F})$, is the maximum number of hyperedges of an  $\mathcal{F}$-free  $r$-uniform hypergraph on $N$ vertices.
 These are  natural generalizations of the classical Tur\'{a}n number for
 graphs \cite{turan}.\\
 
 An $r$-uniform
hypergraph is $\ell$-partite if there is a partition of the vertex set into $\ell$ parts such that each hyperedge
has at most  one vertex in each part. Let $ T_3(N,m-1)$ be the complete $3$-uniform $(m-1)$-partite hypergraph with $N$ vertices and $m-1$  parts $V_1,V_2,...,V_{m-1}$ whose partition sets differ in size by at most $1$. Suppose that  $t_3(N,m-1)$ stands for the number of hyperedges of $T_3(N,m-1)$.\\

 The exact values of $2$-color Ramsey number and Tur\'{a}n number    of $3$-uniform Berge complete hypergraphs  have been  determined by Salia et al. \cite{salia} and Gerbner et al. \cite{gerbcomp}, respectively. Here, we present some needed cases.
 
 \begin{theorem}\em{\cite{salia}}\label{knkm}
 	\begin{equation*}
 	R(\mathcal{B}^3K_m, \mathcal{B}^3K_n)=
 	\left\{
 	\begin{array}{ll}
 	m+n-1 & \ \ \ m=n=3,         \\
 	m+n-2 & \ \ \ m\geq 4,\ n=3\ {\rm or}\ m=n=4,  \\
 	m+n-3 & \ \ \  m\geq 5,\ n=4.
 	\end{array}
 	\right.
 	\quad
 	\end{equation*}
 \end{theorem}

 \begin{theorem}\em{\cite{gerbcomp}}\label{exkm}
	\begin{equation*}
ex(N, \mathcal{B}^3K_m)=
	\left\{
	\begin{array}{ll}
	 5 & \ \ \ m=4,\ N=5,  \\
	t_3(N,m-1) & \ \ \ m=5, N\geq 6.         \\
	\end{array}
	\right.
	\quad
	\end{equation*}
\end{theorem}

They also stated that $T_3(6,4)$ is the largest $\mathcal{B}^3K_5$-free $3$-uniform hypergraph. \\

By Theorem \ref{knkm}, we have the following remark.
\begin{remark}
	$	R(\mathcal{B}^3K_3, \mathcal{B}^3C_3)=5$ and for $m\geq 4$,  $	R(\mathcal{B}^3K_m, \mathcal{B}^3C_3)=m+1$.
\end{remark}


In the following we determine  $R(\mathcal{B}^3K_m, \mathcal{B}^3C_n)$ for $n\in \{4,5\}$ and $m\geq n$.

 \begin{theorem}\label{ck-small}
 	Suppose that  $n \in \{4,5\}$ and $m \geq n$.  Then $R(\mathcal{B}^3K_m, \mathcal{B}^3C_n)= m+1$.
 \end{theorem}
\noindent\textbf{proof. }

To see that $R(\mathcal{B}^3K_m, \mathcal{B}^3C_n)\geq m+1$,  $n \in \{4,5\}$ and $m\geq n$, consider two sets $A$ and $B$, where $|A|=2$ and $|B|=m-2.$ Now,  color all triples crossing  $B$ in at least two points by red and the remaining triples by blue. Obviously, there is no red  $\mathcal{K}_{m}^{(3)}$ and no blue $\mathcal{C}_{n}^{(3)}$, $n\in\{4,5\}$, in this coloring.\\
Since $ex(5,\mathcal{B}^3K_4)=5$ ( Theorem \ref{knkm}) and $ex(5,\mathcal{B}^3C_4)\leq 4$ ( Remark \ref{k5-6e}), we have $R(\mathcal{B}^3K_4, \mathcal{B}^3C_4)\leq 5$. Also, for  $m\geq 5$, using Theorem \ref{knkm}, $R(\mathcal{B}^3K_m, \mathcal{B}^3C_4)\leq m+1$.
 Hence,  $R(\mathcal{B}^3K_m, \mathcal{B}^3C_4)=m+1$ for $m\geq 4$.
	Now, we show that for $m\geq 5$,  $R(\mathcal{B}^3K_m, \mathcal{B}^3C_5)\leq m+1$.\\
	 Suppose that   $\mathcal{H}=\mathcal{K}_{m+1}^3$ is $2$-hyperedge  colored with colors red and blue. It will be proven that $\mathcal{K}_m^{(3)}\subseteq \mathcal{H}_{\textbf{r}}$ or $\mathcal{C}_5^{(3)}\subseteq \mathcal{H}_{\textbf{b}}$.  
	\begin{itemize}
\item[(i)] $m=5$.\\
 Let $V(\mathcal{H})=\{v_1,v_2,\ldots, v_6\}$.  We may assume that $ \mathcal{H}_\textbf{r}$ is   $\mathcal{B}^3K_5$-free.
By Theorem \ref{exkm}, we have  $ex(6,\mathcal{B}^3K_5)=t_3(6,4)=12$. First, let $|E(\mathcal{H}_{\textbf{r}})|=12$ and  $T_1=\{v_1,v_2\}$, $T_2=\{v_3,v_4\}$, $T_4=\{v_5\}$ and $T_4=\{v_6\}$ be parts of   $T_3(6,4)$.  Therefore,  $\mathcal{H}_{\textbf{b}}$ contains a  $\mathcal{C}_{5}^{(3)}$ with the core sequence $v_1, v_3, v_5,v_2,v_4$ and hyperedges  $\{v_1,v_2,v_3\}$, $\{v_3,v_4,v_5\}$, $\{v_1,v_2,v_5\}$, $\{v_1,v_2,v_4\}$,  $\{v_1,v_3,v_4\}$. That holds the assertion.\\
 Now, let  $|E(\mathcal{H}_{\textbf{r}})|<12$. Thus, $|E(\mathcal{H}_{\textbf{b}})|\geq 9$. Since $R(\mathcal{B}^3K_5, \mathcal{B}^3C_4)=6$, we can assume that   $\mathcal{C}_4^{(3)}\subseteq \mathcal{H}_{\textbf{b}}$. Let $v_1,v_2,v_3,v_4$ be the core sequence of   $\mathcal{C}=\mathcal{C}_4^{(3)}$. If $|E(\mathcal{H}_{\textbf{b}}[\{v_1,v_2,v_3,v_4,v_5\}])|\geq 7$ or $|E(\mathcal{H}_{\textbf{b}}[\{v_1,v_2,v_3,v_4,v_6\}])|\geq 7$, By Lemma \ref{k5-3e}, we are done. Otherwise,  one can easily check that $d_{\mathcal{H}_{\textbf{b}}}(v_5)\geq 3$ and $d_{\mathcal{H}_{\textbf{b}}}(v_6)\geq 3$. By some simple discussions on different cases it can be shown that the cycle $\mathcal{C}$ can be extended to a  $\mathcal{C}_5^{(3)}$ in $\mathcal{H}_{\textbf{b}}$. That completes the proof.
 
 \item[(ii)]  $m\geq 6$.\\
 In this case, we use induction on $m$. 
 Since  $R(\mathcal{B}^3K_{m-1}, \mathcal{B}^3C_5)=m$, we can assume that   $\mathcal{K}_{m-1}^{(3)}\subseteq \mathcal{H}_{\textbf{r}}- v$ for every $v\in V(\mathcal{H})$.\\
  First, assume that there is a vertex $u\in V(\mathcal{H})$ such that $\textbf{r}\in c(uv)$ in $\Gamma_3(\mathcal{H})$, for every $v\in V(\mathcal{H})\setminus \{u\}$. Furthermore,  assume that $\mathcal{K}=  \mathcal{K}_{m-1}^3 \subseteq \mathcal{H}_{\textbf{r}}-u$.  Let $\{v_1,v_2,\ldots, v_{m-1}\}$ and $\{h_{ij} \ :\ 1\leq i<j\leq m-1\}$ be the set of main vertices and hyperedges of $\mathcal{K}$, respectively.
Set  $A_i=\{h\ :\ h\in E(\mathcal{H}_{\textbf{r}}), \{u,v_i \}\subseteq h\}$, for $1\leq i\leq m-1$.
Now, let $G=[X,Y]$ be the bipartite graph defined as follows. Let $X=\{uv_i \ :\ 1\leq i\leq m-1\}$ and $Y=\cup_{i=1}^{m-1}A_i$. For every $uv_i \in X$ and $h \in Y$, $uv_i\sim_{G}h$,  if and only if $h\in A_i$. Clearly,   $d_{G}(uv_i) \geq 3$ for every $1\leq i \leq m-1$. Also, we have   $d_{G}(h) \leq 2$ for every $h \in Y$. Using the Hall's Theorem,  $G$ contains the matching $M$ of $X$. Let $Y'\subseteq Y$ be the members of $Y$ saturated by $M$. Clearly,  $\{u,v_1,v_2,\ldots, v_{m-1}\}$ and $\{h_{ij}\ :\ 1\leq i<j\leq m-1\}\cup Y'$ makes     a complete Berge hypergraph $\mathcal{K}_m^{(3)}$ in $\mathcal{H}_{\textbf{r}}$.\\
Therefore, for every vertex $v\in V(\mathcal{H})$, there is a vertex $v'\in V(\mathcal{H})\setminus \{v\}$ such that $c(vv')=\{\textbf{b}\}$ in $\Gamma_3(\mathcal{H})$.  Since  $R(\mathcal{B}^3K_{m}, \mathcal{B}^3C_4)=m+1$, we may assume that $\mathcal{C}_4^{(3)}\subseteq \mathcal{H}_{\textbf{b}}$. Otherwise,  $\mathcal{K}_m^{(3)}\subseteq \mathcal{H}_{\textbf{r}}$ and the proof is completed.  Set  $\mathcal{C}=\mathcal{C}_4^{(3)}$ and  let $v_1,v_2,v_3,v_4$ and $h_1,h_2,h_3,h_4$ be the core sequence and hyperedges of   $\mathcal{C}$, respectively. Also, let $U=V(\mathcal{H})\setminus \{v_1,v_2,v_3,v_4\}$. Clearly, $|U|\geq 3$.  It suffices to consider the following three cases to prove that  $\mathcal{C}_5^{(3)}\subseteq \mathcal{H}_{\textbf{b}}$.
\\

(a) $c(uu')=\{\textbf{b}\}$ in $\Gamma_3(\mathcal{H})$, for some $u,u'\in U$. \\
If for some $1\leq t \leq 4$, there are $v_t$ and $v_{t+1}$ such that $\{u,u',v_t\}$ and $\{u,u',v_{t+1}\}$ are blue, then we have a blue $\mathcal{C}_5^{(3)}$ by embedding $u$ between $v_t$ and $v_{t+1}$, removing $h_t$ and adding $\{u,u',v_t\}$ and $\{u,u',v_{t+1}\}$. Otherwise, for some $1\leq t \leq 4$, hyperedges $\{u,u',v_t\}$ and $\{u,u',v_{t+2}\}$ are blue. Since $c(uu')=\{\textbf{b}\}$, there exists a third hyperedge in $  E(\mathcal{H}_{\textbf{b}})$, say $h'$, containing $u$ and $u'$. Therefore, $v_{t-1},v_t,u,u',v_{t+2}$ represents the core sequence of a blue $\mathcal{C}_5^{(3)}$ with hyperedges  $h_{t-1}, 
\{u,u',v_t\}, h', \{u,u',v_{t+2}\},h_{t+2}$.\\

(b) $c(uv_{\ell})=c(u'v_{\ell})=\{\textbf{b}\}$ in $\Gamma_3(\mathcal{H})$, for some $u,u'\in U$ and $1\leq \ell \leq 4$. \\
If $\{u,v_{\ell},v_{\ell+1}\}$ is blue, we have a blue $\mathcal{C}_5^{(3)}$ with the core sequence $v_{\ell-1},v_{\ell},u,v_{\ell+1},v_{\ell+2}$ and the hyperedges  $h_{\ell-1}, 
 h',\{u,v_{\ell},v_{\ell+1}\},h_{\ell+1},h_{\ell+2}$, where $h' \in  E(\mathcal{H}_{\textbf{b}})$ is a  hyperedge containing $u$ and $v_{\ell}$ distinct from $h_{\ell -1}$ and $\{u,v_{\ell},v_{\ell+1}\}$. Therefore, we may assume that $\{u,v_{\ell},v_{\ell+1}\}$, $\{u,v_{\ell-1},v_{\ell}\}$, $\{u',v_{\ell},v_{\ell+1}\}$ and $\{u',v_{\ell-1},v_{\ell}\}$ are red. Hence, $\{u,u',v_{\ell}\}$, $\{u,v_{\ell},v_{\ell+2}\}$ and $\{u',v_{\ell},v_{\ell+2}\}$ are blue by the assumption  $c(uv_{\ell})=c(u'v_{\ell})=\{\textbf{b}\}$. Therefore, $\mathcal{C}$ can be extended to a  $\mathcal{C}_5^{(3)}$ in $\mathcal{H}_{\textbf{b}}$ with the core sequence  $v_{\ell},v_{\ell+1},v_{\ell+2},u,u'$ and the hyperedges  $h_{\ell},h_{\ell+1},\{u,v_{\ell},v_{\ell+2}\},\{u,u',v_{\ell}\},\{u',v_{\ell},v_{\ell+2}\}$.\\
 
(c) $c(uv_{\ell})=c(u'v_{\ell+1})=\{\textbf{b}\}$ in $\Gamma_3(\mathcal{H})$, for some $u,u'\in U$ and $1\leq \ell \leq 4$. \\
 By an argument similar to the previous case, we can assume that  $\{u,v_{\ell},v_{\ell+1}\}$, $\{u,v_{\ell-1},v_{\ell}\}$, $\{u',v_{\ell},v_{\ell+1}\}$ and $\{u',v_{\ell+1},v_{\ell+2}\}$ are red. Thus, $\{u,u',v_{\ell+1}\}$ is blue.  In this case,  sequence $v_{\ell-1},v_{\ell},u,v_{\ell+1}, v_{\ell+2}$ and hyperedges $h_{\ell-1},h',\{u,u',v_{\ell+1}\},h_{\ell+1},h_{\ell+2}$ make a blue $\mathcal{C}_5^{(3)}$, where $h'$ is a blue hyperedge containing $u$ and $v_{\ell}$. Hence, the assertion holds.
$\hfill\blacksquare$\\
\end{itemize}

In the sequel,  we will present the following lemma  and then use it to determine $R(\mathcal{B}^3K_m, \mathcal{B}^3C_n)$ for $m\geq n\geq 6$ and $m  \geq 11$.
\begin{lemma}\label{exten2}
	Consider a $t$-hyperedge coloring of   $3$-uniform hypergraph $\mathcal{H}$  with colors $\{\mathfrak{c}_1,\ldots,\mathfrak{c}_t\}$.
	If $K_m \subseteq \Gamma_3(\mathcal{H})$ and  there is a color $\mathfrak{c}_{\ell}$ for some  $1\leq \ell \leq t$ such that for every $e \in E(K_m)$,  $\mathfrak{c}_{\ell}\in c(e)$, then $\mathcal{K}_m^{(3)} \subseteq  \mathcal{H}_{\mathfrak{c}_{\ell}}$.
	%
\end{lemma}	
\noindent\textbf{proof. }Suppose that $V(K_m)=\{v_1,v_2,\ldots, v_m\}$.  Let $G^*=[X,Y]$ be the bipartite graph defined as follows. Let $X=\{v_{ij}\ :\ 1\leq i< j \leq m\}$ and $Y=\{h\ :\ h\in E(\mathcal{H}_{\mathfrak{c}_{\ell}}) \}$. For every $v_{ij} \in X$ and $h \in Y$, $v_{ij}\sim_{G^*}h$,  if and only if $\{v_i,v_j\}\subseteq h$. Clearly,   $d_{G^*}(v_{ij}) \geq 3$ for every $1\leq i<j \leq m$. Also, we have   $d_{G^*}(h) \leq 3$ for every $h \in Y$. Using the Hall's theorem,  $G^*$ contains a matching of $X$ that leads to  a complete Berge hypergraph $\mathcal{K}_m^{(3)}$ in $\mathcal{H}_{\mathfrak{c}_{\ell}}$.
$\hfill\blacksquare$

 \begin{theorem}\label{ck-general}
 	Let $m\geq n\geq 6$ and $m  \geq 11$. Then $R(\mathcal{B}^3K_m, \mathcal{B}^3C_n)= m+\lfloor \frac{n-1}{2}\rfloor -1.$
 \end{theorem}
\noindent\textbf{proof. }First we show that $R(\mathcal{B}^3K_m, \mathcal{B}^3C_n)\geq m+\lfloor \frac{n-1}{2}\rfloor -1.$ 
	   To see that, consider  two sets $A$ and $B$, where $|A|=\lfloor \frac{n-1}{2}\rfloor$ and $|B|=m-2.$ Now,  color all triples intersecting  $B$ in at least two points by red and the other triples by blue.  Clearly,  there is neither red $\mathcal{K}_m^{(3)}$ nor blue $\mathcal{C}_n^{(3)}$ in this coloring.
	Therefore, it suffices   to  show that $R(\mathcal{B}^3K_m, \mathcal{B}^3C_n)\leq m+\lfloor \frac{n-1}{2}\rfloor -1.$\\
	 Now, suppose that  $\mathcal{H}=\mathcal{K}_{m+\lfloor \frac{n-1}{2}\rfloor -1}^3$ is $2$-hyperedge colored  with colors red and blue.  We show that $\mathcal{K}_m^{(3)}\subseteq \mathcal{H}_{\textbf{r}}$ or $\mathcal{C}_n^{(3)}\subseteq \mathcal{H}_{\textbf{b}}$.
	Using Theorem \ref{ck-small}, we have  $\mathcal{K}_{m}^{(3)} \subseteq \mathcal{H}_{\textbf{r}}$ or   $\mathcal{C}_{5}^{(3)} \subseteq \mathcal{H}_{\textbf{b}}$. If $\mathcal{K}_{m}^{(3)} \subseteq \mathcal{H}_{\textbf{r}}$, the proof is completed. Otherwise, let $\mathcal{C}_{\ell}^{(3)}$ be a Berge cycle with the maximum length in $\mathcal{H}_{\textbf{b}}$. Also, let $v_1,v_2,\ldots, v_{\ell}$ be the   core sequence of $\mathcal{C}_{\ell}^{(3)}$, where $5\leq \ell \leq n$.
	 If $\ell =n$, we are done. Otherwise, we may assume that $5 \leq  \ell \leq n-1$. Set $\mathcal{C}=\mathcal{C}_{\ell}^{(3)}$ and $U=V(\mathcal{H})\setminus \{v_1,v_2,\ldots,v_{\ell}\}$.
	 Now, we have the following claim.
	 
\begin{claim}\label{e}
	 For every $u\in U$ and   $e=uv\in E(\Gamma_3(\mathcal{H}))$ we have  $c(e)\neq \{{\bf b}\}$.	
\end{claim}
\noindent\textbf{proof of Claim \ref{e}. }Suppose to the contrary that there is an edge  $e=uv\in E(\Gamma_3(\mathcal{H}))$ such that $u\in U$  and $c(e)= \{{\bf b}\}$. First, let $v\in \{v_1,\ldots,v_{\ell}\}$. So, $v=v_i$ for some $1\leq i \leq \ell$. Also, let $h_1,h_2,\ldots,h_{\ell}$ be the hyperedges of $\mathcal{C}$ such that every $h_i$ contains $v_i$ and $v_{i+1}$. Assume that the  hyperedge $f_1=\{u,v_{i-1},v_{i}\}$  is blue. Since there is a blue hyperedge $f_2\neq f_1$ containing $u$ and $v_i$,   we have a  $\mathcal{C}_{\ell+1}^{(3)}\subseteq \mathcal{H}_\textbf{b}$, by setting $u$ between  $v_{i-1}$ and $v_{i}$, removing hyperedge $h_{i-1}$  and adding  hyperedges $f_1$ and  $f_2$, a contradiction. Therefore, we may assume that $\{u,v_{i-1},v_{i}\}$  and $\{u,v_i,v_{i+1}\}$ are red.\\
	Since $\ell \geq 5$ and $c(e)= \{{\bf b}\}$, there exist two consecutive vertices $v_t$ and $v_{t+1}$ on $\mathcal{C}$, $t \notin \{i-2,i-1,i,i+1\}$ such that $\{u,v_i,v_{t}\}$ and $\{u,v_i,v_{t+1}\}$ are blue. Clearly, we have a blue  $\mathcal{C}_{\ell+1}^{(3)}$, by embedding $u$ between $v_t$ and $v_{t+1}$, removing $h_t$ and using hyperedges $\{u,v_i,v_{t}\}$ and $\{u,v_i,v_{t+1}\}$. It contradicts the assumption that $\mathcal{C}$ is maximum.\\ 
	Now, let  $v\in U$. Since $c(e)=\{\textbf{b}\}$, there are at most two vertices $v_j$ and $v_{j'}$ on $\mathcal{C}$ such that hyperedges $\{u,v,v_{j}\}$ and $\{u,v,v_{j'}\}$ are red. The fact $\ell \geq 5$ implies that there are two consecutive vertices $v_t$ and $v_{t+1}$ on $\mathcal{C}$  where $\{u,v,v_{t}\}$ and $\{u,v,v_{t+1}\}$ are blue. One can  find a  $\mathcal{C}_{\ell+1}^{(3)}$ in $\mathcal{H}_\textbf{b}$   by an argument  similar to the above, a contradiction.
$\hfill\square$\\

By the Claim \ref{e}, for every  edge $e=uv \in \Gamma_3(\mathcal{H})$ such that  $|\{u,v\}\cap U| \geq 1$, we have   $\textbf{r} \in c(e)$. 
	Let $e$ and $e'$ be two edges of $\Gamma_3(\mathcal{H})$ such that $c(e)=c(e')=\{\textbf{b}\}$. First, let $e$ and $e'$ are disjoint. Set $e=v_iv_{j}$ and $e'=v_{i'}v_{j'}$, 
	$1 \leq i,i',j,j' \leq \ell$.     Since $m \geq 11$ and $\ell \leq n-1$, we have  $|U| \geq 5$. Therefore, there is a vertex $w \in U$ such that the   hyperedges $\{w, v_i, v_{j}\}$ and $\{w,v_{i'}, v_{j'}\}$ are blue. If one of the set $\{v_i,v_{i'}\}$, $\{v_i,v_{j'}\}$, $\{v_{j},v_{i'}\}$ or $\{v_{j},v_{j'}\}$ is placed consecutively on the core sequence of $\mathcal{C}$, for instance $v_i=v_t$ and $v_{i'}=v_{t+1}$ for some $1\leq t \leq \ell $,  a  $\mathcal{C}_{\ell +1}^{(3)}$ is obtained in $\mathcal{H}_{\textbf{b}}$  by embedding $w$ between $v_t$ and $v_{t+1}$. It contradicts the assumption that $\mathcal{C}$ is maximum.\\
	Now, let $e=v_iv_j$ and $e'=v_iv_{j'}$.  Since  $|U| \geq 5$, similarly, one can see that there is a vertex $w' \in U$ such that the   hyperedges $\{w', v_i, v_{j}\}$ and $\{w',v_{i}, v_{j'}\}$ are blue. If $v_j$ and $v_{j'}$ are placed consecutively  on the core sequence of $\mathcal{C}$, we can embed $w'$ between $v_j$ and $v_{j'}$ and use hyperedges  $\{w', v_i, v_{j}\}$ and $\{w',v_{i}, v_{j'}\}$ to construct a blue $\mathcal{C}_{\ell +1}^{(3)}$, a contradiction.\\
	Let  $\{g_1, g_2, \ldots, g_k\}\subseteq E(\Gamma_3(\mathcal{H}))$ be the set of all edges  such that $c(g_i)=\{\textbf{b}\}$, for every $1\leq i \leq k$. Also, let $G$ be the subgraph of $\Gamma_3(\mathcal{H})$ induced by $\{g_1, g_2, \ldots, g_k\}$. By the above discussion,   $|V(G)|\leq  \lfloor \frac{\ell}{2} \rfloor \leq \lfloor\frac{n-1}{2}\rfloor$.\\
	 Now, set $Y= U \cup (\{v_1,\ldots,v_{\ell}\}\setminus V(G))\cup \{z\}$ for some $z \in V(G)$. One can easily see that $|Y|\geq m$. Let $\{y_1,y_2,\ldots, y_m\} \subseteq Y.$ Since for every $1\leq i,j \leq m$,    $\textbf{r} \in c(y_iy_j)$ in $\Gamma_3(\mathcal{H})$,  by the Lemma \ref{exten2},  we have  $\mathcal{K}_m^{(3)}\subseteq \mathcal{H}_{\textbf{r}}$. That completes the proof.
	$\hfill\blacksquare$\\


\begin{thebibliography}{10}
	
	\bibitem{axenovich}
	M. Axenovich and A.  Gy\'{a}rf\'{a}s, \textit{A note on Ramsey numbers for Berge-G hypergraphs},  Discrete Math., 342 (2019), pp. 1245-1252.
	
	\bibitem{davoodi}
	A. Davoodi, E. Gy\H{o}ri, A. Methuku and C. Tompkins, \textit{An Erd\H{o}s-Gallai type theorem for uniform hypergraphs},  European J. Combin., 69 (2018), pp. 59-162.
	
	\bibitem{gerbram}
	D. Gerbner, \textit{On Berge-Ramsey problems}, Electron. J. Combin., 27 (2020), P2.39. 
	
	\bibitem{omidi}
	D. Gerbner, A. Methuku, G. Omidi, and M. Vizer, \textit{Ramsey problems for Berge
	hypergraphs},   SIAM J. Discrete Math., 34 (2020), pp. 351-369.
	
	\bibitem{gerbcomp}
	D. Gerbner, A. Methuku and C. Palmer, \textit{General lemmas for Berge-Tur\'{a}n hypergraph problems},  European J. Combin., 68 (2020), 103082.
	
	\bibitem{vizer}
	D. Gerbner, A. Methuku and M. Vizer,   \textit{Asymptotics for the Tur\'{a}n number of Berge-$K_{2,t}$}, J. Combin.
		Theory Ser. B, 137 (2019), pp.  264-290.
	
	\bibitem{gerbdef}
	D. Gerbner and C. Palmer, \textit{Extremal Results for Berge Hypergraphs}, SIAM J.
		Discrete Math., 31 (2017), pp.  2314-2327.
	
	\bibitem{gyarfcomp}
	A. Gy\'{a}rf\'{a}s, \textit{The Tur\'{a}n number of Berge-$K_4$ in triple systems}, SIAM J. Discrete Math., 33 (2019), pp.  383-392.
	
	\bibitem{gyarfconj1}
	A. Gy\'{a}rf\'{a}s, J. Lehel, G. N. S\'{a}rk\"{o}zy, and R. H. Schelp, \textit{Monochromatic Hamiltonian Berge-cycles in colored complete uniform hypergraphs},  J. Combin.
		Theory Ser. B, 98 (2008), pp. 342-358.
	
	\bibitem{gyarfconj2}
	A. Gy\'{a}rf\'{a}s and G. N.  S\'{a}rk\"{o}zy, \textit{The 3-colour Ramsey number of a 3-uniform
	Berge cycle},  Combin. Probab. Comput., 20 (2011), pp. 53-71.
	
	\bibitem{zemeredi}
	A. Gy\'{a}rf\'{a}s, G.N. S\'{a}rk\"{o}zy and E. Szemer\'{e}di, \textit{Monochromatic matchings in the shadow graph of almost complete hypergraphs}, Ann. Comb., 14 (2010), pp. 245-249.
	
	\bibitem{longcycle} 
	A. Gy\'{a}rf\'{a}s, G.N. S\'{a}rk\"{o}zy and E. Szemer\'{e}di, \textit{Long monochromatic Berge-cycles in colored 4-uniform
	hypergraphs}, Graphs Combin., 26 (2010), pp. 71-76.
	
	\bibitem{gyorik3}
	E. Gy\H{o}ri, \textit{Triangle-free hypergraphs},  Combin. Probab. Comput., 15 (2006), pp. 185-191.
	
	\bibitem{katona}
	E. Gy\H{o}ri, G. Y. Katona and N. Lemons, \textit{Hypergraph extensions of the Erd\H{o}s-Gallai Theorem},
	European J. Combin., 58 (2016), pp. 238-246.
	
	
	
	\bibitem{gyori}
	E. Gy\H{o}ri and N. Lemons, \textit{Hypergraphs with no cycle of a given length}, Combin. Probab. Comput., 21 (2012), pp. 193-201.
	
	\bibitem{mahcycle}
	L. Maherani and G. R. Omidi,  \textit{Monochromatic Hamiltonian Berge-cycles in colored hypergraphs}, Discrete Math., 340 (2017), pp. 2043-2052.
	
	
	\bibitem{mahcomp}
	L. Maherani and M. Shahsiah, \textit{Tur\'{a}n numbers of complete 3-uniform Berge-hypergraphs}, Graphs Combin., 34 (2018), pp.  619-632.
	
	\bibitem{omidicycle}
	G. R. Omidi, \textit{A proof for a conjecture of Gy\'{a}rf\'{a}s, Lehel, S\'{a}rk\"{o}zy and Schelp on Berge-cycles},  Combin. Probab. Comput., 30 (2021), pp. 654-669.
	
	\bibitem{palmer2}
	C. Palmer, M. Tait, C. Timmons and A. Z. Wagner, \textit{Tur\'{a}n numbers for Berge-hypergraphs
	and related extremal problems}, Discrete Math., 342 (2019), pp. 1553-1563. 
	
	
	
	
	
	\bibitem{palvolgy}
	D. P\'{a}lv\"{o}lgyi, \textit{Exponential lower bound for Berge-Ramsey problems},  Graphs Combin., 37 (2021), pp. 1433-1435. 
	
	\bibitem{salia}
	N. Salia, C. Tompkins, Z. Wang and O. Zamora, \textit{Ramsey numbers of Berge hypergraphs and related structures}, Electron. J. Combin. 26 (2019),  P4.40.
	
	\bibitem{turan}
	 P. Tur{\'a}n, \textit{Eine {E}xtremalaufgabe aus der {G}raphentheorie}, Mat. Fiz. Lapok,  48 (1941), pp. 436-452.
	
	
	
\end{thebibliography}
\end{document}